\documentclass[english,10pt]{article}

\usepackage{xcolor}
\usepackage{amsmath}
\usepackage{graphics}
\usepackage{epsfig}
\usepackage{vmargin}
\usepackage{eqnarray}
\usepackage{amssymb}
\usepackage[latin1]{inputenc}
\usepackage{babel,indentfirst,amsfonts,array}

\newtheorem{thm}{Theorem}[section]
\newtheorem{lemme}[thm]{Lemma}
\newtheorem{prop}[thm]{Proposition}

\newtheorem{defi}[thm]{Definition}

\newtheorem{hyp}[thm]{Hypothesis}

\def\P{\mathbb P}
\def\E{\mathbb E}

\def\R{\mathbb R}

\def\Z{\mathbb Z}
\def\N{\mathbb N}
\def\({\left(}
\def\){\right)}
\def\[{\left[}
\def\]{\right]}
\def\fin{\hfill\square}

\def\fin{\hfill $\square$}

\title{\textsc{Random walk in a stratified independent random environment}
\author{Julien Br\'emont}
\date{Universit\'e Paris-Est Cr\'eteil,~novembre 2018}
}

\begin{document}

\maketitle

\setcounter{page}{1}

\begin{abstract}
We study Markov chains on a lattice in a codimension-one stratified independent random environment, exploiting results established in \cite{jb2}. First of all the random walk is transient in dimension at least three. Focusing on dimension two, both recurrence and transience can happen, but transience remains by far the most general situation. We identify the critical scale of the local drift along the strata corresponding to the frontier between the two regimes. \end{abstract}

\footnote{
\begin{tabular}{l}\textit{AMS $2000$ subject classifications~: 60G17, 60J10, 60K37.} \\
\textit{Key words and phrases~: Markov chain, recurrence, stratification, random environment, independence.} 
\end{tabular}}

\section{Introduction}
A very first and important question on the asymptotic behaviour of Markov chains on a lattice in heterogeneous environment is the question of recurrence/transience. We focus in this paper on the situation where the transition laws depend on a single coordinate, continuing the line of research developed in \cite{jb1,jb2}. Historically, a PDE model has first been introduced by Matheron and de Marsily \cite{mama} at the beginning of the eighties and a probabilistic version was introduced by Campanino and Petritis \cite{cp} in 2003. Let us detail the extended version studied in \cite{jb2}.

\medskip
\noindent
Fixing $d\geq1$, we consider a Markov chain $(S_n)_{n\geq0}$ in $\Z^{d}\times\Z$, with $S_0=0$. Quantities relative to the first (resp. second) coordinate in $\Z^d$ (resp. $\Z$) are declared ``horizontal" (resp. ``vertical"). We assume that the transition laws are constant on each affine hyperplane $\Z^d\times\{n\}$, $n\in \Z$. To model this, for each vertical $n\in\Z$, let positive reals $p_n,q_n,r_n$, with $p_n+q_n+r_n=1$, and a probability measure $\mu_n$ with support in $\Z^d$, satisfying the following conditions~:

\begin{hyp}
\label{hypo} $\exists\eta>0$, $\forall n\in\Z$, $\min\{p_n,q_n,r_n\}\geq\eta$, $\sum_{k\in\Z^d}\|k\|^{\max(d,3)}\mu_n(k)\leq1/\eta$ and the spectrum of $\sum_{k\in\Z^d}kk^t\mu_n(k)$ is included in $[\eta,+\infty).$
\end{hyp}

\noindent
Observe that when $d=1$, the last condition reduces to $\mu_n(0)\leq 1-\eta$. The local horizontal drift is by definition $\varepsilon_n=\sum_{k\in\Z^d}k\mu_n(k)$, the expectation of $\mu_n$, $n\in\Z$. The transition laws are now defined, for all $(m,n)\in\Z^d\times\Z$ and $ k\in\Z^d$, by~:

$$(m,n)\overset{p_n}{\longrightarrow}(m,n+1),~(m,n)\overset{q_n}{\longrightarrow}(m,n-1),~(m,n)\overset{r_n\mu_n(k)}{\longrightarrow}(m+k,n).$$

\medskip
The Matheron-De Marsily or Campanino-Petritis model corresponds to taking $d=1$, with $p_n=q_n=r_n=1/3$ and $\mu_n=\delta_{\alpha_n}$, fixing some $(\alpha_n)_{n\in\Z}\in\{\pm 1\}^{\Z}$, corresponding to orientating horizontal lines. It is shown in \cite{cp}, for example, that the random walk is recurrent when $\alpha_n=(-1)^n$ and transient when $(\alpha_n)$ are independent and identically distributed (we say $i.i.d$ for the sequel) random variables with $\P(\alpha_n=\pm1)=1/2$, for a.-e. realization. Several extensions around the orientation setting then followed; see Campanino-Petritis \cite{cp2} for a review. In \cite{jb1}, for the model introduced in the above paragraph, with $d=1$ and $p_n=q_n$, $\forall n\in\Z$, (vertically flat model) a characterization of the recurrence/transience of the random walk was given. The rather abstract form of this criterion directly comes from the computation of a Poisson kernel in the half-plane, as appearing in the recurrence criterion for $i.i.d.$ random walks on $\Z^d$ (see Spitzer \cite{spitz}). The latter is an essential ingredient in the proof. The principal result of \cite{jb1} is the main step in analyzing the recurrence/transience properties of the vertically flat model, but some (little) extra work is needed to treat concrete examples. This was done at the end of \cite{jb1}, allowing to prove accurate results about recurrence/transience for this model. 

\medskip
\noindent
The same recurrence criterion as in \cite{jb1} was next extended in \cite{jb2} to the general setting introduced at the beginning of this section. The structure stays the same, but a new metrization of the environment has to introduced, informally corresponding to the way the random walk sees the environment. An interpretation of the criterion was given in relation with the level lines of some interesting function of two variables measuring horizontal dispersion. This will be detailed later. Several examples of recurrent and transient random walks were finally given, both in dimensions 2 and 3 (i.e. $d\in\{1,2\}$). In dimension $\geq4$ (i.e. $d=3$), the random walk appears to be always transient (\cite{jb2}, prop. 2.5, 1)$i)$). 

\medskip
The purpose of the present short article is to specifically study the important case when the stratifications of the environment are random and independent, with a quenched point of view. We show that the recurrence criterion of \cite{jb2}, once mastered the user's manual, gives sharp quantitative results, probably having some flavour of what should happen for much more general models of random walks in independent environments. 

\medskip
\noindent
Let us present the results. We always assume hypothesis \ref{hypo}. In the sequel, randomness is for the environment; we never enter the mechanism of the random walk. 

\begin{prop} 
\label{prop1} 

$ $

\noindent
If$(p_n/q_n)_{n\in\Z}$ are i.i.d. random variables, with either $\E(\log (p_0/q_0))\not=0$ or $d\geq 2$, then for a.-e. realization the random walk is transient.
\end{prop}
 
\noindent
This will follow rather readily from \cite{jb2}. Indeed, a key observation, already made in \cite{cp} and direct consequence of the invariance of the environment under horizontal translations, is that the random walk restricted to vertical jumps is a Markov chain on $\Z$ (this is never true in general, for a non-stratified environment). Its transition laws are~:

$$n\overset{p_n/(p_n+q_n)}{\longrightarrow}n+1,~n\overset{q_n/(p_n+q_n)}{\longrightarrow}n-1.$$ 

\medskip
\noindent
We call it the ``vertical random walk''. Classically when $\E(\log (p_0/q_0))\not=0$, the vertical random walk is transient (see \cite{sol}), so is the initial random walk. Transience when $\E(\log (p_0/q_0))=0$ and $d=2$ will also be shortly proved. Recall here in passing the related conjecture that any random walk in $i.i.d.$ random environment on $\Z^3$, with ellipticity conditions on the data, is transient.

\medskip
To complete the study of the independent case, we now focus on the planar case, i.e. $d=1$, with $\E(\log (p_0/q_0))=0$. Recall that $\varepsilon_n=\sum_{k\in\Z}k\mu_n(k)$.

\medskip
\noindent
When $\P(p_0=q_0)=1$, i.e. simply $p_n=q_n$, $n\in\Z$, as a consequence of \cite{jb1} prop. 1.4 $i)$, the random walk is recurrent for any $(\varepsilon_n)_{n\in\Z}$ verifying~:

$$|\sum_{-n\leq k\leq 0}\varepsilon_k|+|\sum_{0\leq k\leq n}\varepsilon_k|=O((\log n)^{1/2}).$$

\noindent
In the other direction, cf \cite{jb1} prop. 1.6, when the $(\varepsilon_n)_n$ are independent random variables such that for some $\varepsilon>0$~:

$$\liminf \frac{1}{N}\sum_{n=1}^N\P\({|\sum_{0\leq k\leq n}\varepsilon_k|\geq (\log n)^{1+\varepsilon}}\)>0,$$

\noindent
then, for a.-e. realization, the random walk is transient. This extended former results of \cite{cp} and Devulder-P\`ene \cite{DP}. In the independent case, the critical size of the sums $\sum_{0\leq k\leq n}\varepsilon_k$ with respect to recurrence/transience is expected to be $\log n$, so corresponding to $\varepsilon_n$ having order $1/n$.

\medskip
We now consider the remaining case, reformulating hypothesis \ref{hypo} in a slightly stronger form. We denote by $B(x,r)$ the interval $]x-r,x+r[\subset\R$.
 
\begin{thm}
\label{indep} 

$ $

\noindent
Let $d=1$ and assume $\exists\eta>0,\forall n\in\Z$, $\min\{p_n,q_n,r_n\}\geq\eta$, $\mbox{Supp}(\mu_n)\subset[-1/\eta,1/\eta]$, $\mu_n(0)\leq 1-\eta$, with $(p_n/q_n)_{n\in\Z}$ non-constant $i.i.d.$ random variables, $\E(\log (p_0/q_0))=0$.

\medskip
\noindent 
i) For a.-e. realization~: any $(\varepsilon_n)_{n\in\Z}$ such that $\varepsilon_n=O(\exp(-|n|^{1/2+\delta}))$, $\delta>0$, imply that the random walk is recurrent; any $(\varepsilon_n)_{n\in\Z}$ such that for some $\delta>0$ and large $n$, $\varepsilon_n\geq \exp(-|n|^{1/2-\delta})$ (hence $>0$) imply that the random walk is transient. 

\medskip
\noindent ii) Suppose that the $(\varepsilon_n)_n$ are independent random variables, independent of the $(p_n/q_n)_{n\in\Z}$. Assume $\exists\delta>0$ so that for large $n\in\N$ and all $y\in\R$, $\P(\varepsilon_n\not\in B(y,\exp(-n^{1/2-\delta})))\geq n^{-\delta/5}$. Then for a.-e. realization the random walk is transient.

\end{thm}

\bigskip
\noindent 
\begin{remark}
We shall prove versions a little stronger, as indicated at the beginning of the proofs below. Roughly when $\P(p_0=q_0)<1$, the critical scale for $\epsilon_n$ with respect to recurrence/transience is $\exp(-|n|^{1/2})$, hence much smaller than for the vertically flat case. The vertical random walk is now Sina\"\i's random walk, not simple random walk. It stays a long time in deep valleys, coming back to zero with very small frequency. For the original random walk, the landscape looks like a succession of horizontally invariant canyons, where it stays confined for a long time. Just a little of horizontal drift is enough to make the random walk transient.

\end{remark}

\bigskip
\noindent 
\begin{remark} Item $ii)$, when $\P(\mu_n=\delta_1)=\P(\mu_n=\delta_{-1})=1/2$, giving $\P(\varepsilon_n=\pm1)=1/2$, was first proved by Kochler in his Ph-D thesis \cite{koch}, after a long and delicate analysis of the Brownian path. 
\end{remark}

\section{Preliminaries}

Let us begin with some notations. We first fix an integer $K>2(1+1/\eta)$.

\begin{defi} 

$ $

\noindent
i) Set $a_n=q_n/p_n$ and~:

$$\rho_n=\left\{{\begin{array}{cc}
a_1\cdots a_n,&n\geq1,\\
1,&n=0,\\
a_{n+1}\cdots a_{-1}a_0&~n\leq-1.\end{array}
}\right.$$

\smallskip
\noindent
For $n\geq0$, let~:

$$v_+(n)=\sum_{0\leq k\leq n}\rho_k\mbox{ and }v_-(n)=a_0\sum_{-n-1\leq k\leq -1}\rho_k.$$

\smallskip
\noindent
ii) For functions $f(n)$ and $g(n)$ defined on $\N=\{0,1,\cdots\}$, we write $f\asymp g$ if there exists a constant $C>0$ so that for large $n$, $(1/C)f(n)\leq g(n)\leq Cf(n)$. We write $g\preceq f$ if $g(n)\leq Cf(n)$ only.

\smallskip
\noindent
iii) Let $f:\N\rightarrow\R_+$, increasing, with $\lim f(n)=+\infty$. For large enough $x\in\R_+$, let $f^{-1}(x)=\max\{n\in\N~|~f(n)\leq x\}$. Note that $f(f^{-1}(x))\leq x<f(f^{-1}(x)+1)$ and $f^{-1}(f(n))=n$.

\end{defi}

Notice that $\rho_{n+1}/\rho_{n}\in[\eta,1/\eta]$, because $a_k\in[\eta,1/\eta]$, $\forall k$. It is immediate to see that, a.-e., $v_+(n)$ and $v_-(n)$ are both increasing and tend to $+\infty$, when $\E(\log (p_0/q_0))=0$. For example, the $i.i.d.$ random walk $\sum_{0\leq k\leq n}\log (q_k/p_k)$ has an integrable and centered step, hence is recurrent, so almost-surely $\rho_n$ does not go to $0$. 

\medskip
Let us now introduce the functions describing the average horizontal macrodispersion of the environment.

\begin{defi} 

$ $

\noindent
i) The structure function, depending only on the vertical, is defined for $n\geq0$ by~:

$$\Phi_{str}(n)=\({n\sum_{-v_-^{-1}(n)\leq k\leq v_+^{-1}(n)}1/\rho_k}\)^{1/2}.$$

\medskip
\noindent
2) For $d=1$ and $m,n\geq0$, introduce~:
$$\Phi(-m,n)=\({\sum_{-v_-^{-1}(m)\leq k\leq l\leq v_+^{-1}(n)}\rho_k\rho_l\[{1/\rho^2_k+1/\rho_l^2+\({\sum_{s=k}^l\frac{r_s\varepsilon_s}{p_s\rho_s}}\)^2}\]}\)^{1/2}.$$

\noindent
For $n\geq 0$, set $\Phi(n)=\Phi(-n,n)$ and $\Phi_{+}(n)=\sqrt{\Phi^2(-n,0)+\Phi^2(0,n)}$. 
\end{defi}

Observe that for $n\geq1$ :

\begin{eqnarray}
\sum_{0\leq k\leq l\leq v_+^{-1}(n)}\rho_k\rho_l\({1/\rho^2_k+1/\rho_l^2}\)&=&\sum_{0\leq k\leq l\leq v_+^{-1}(n)}(\rho_k/\rho_l+\rho_l/\rho_k)\nonumber\\
&\asymp&\sum_{0\leq k\leq v_+^{-1}(n)}\rho_k\sum_{0\leq l\leq v_+^{-1}(n)}1/\rho_l.\nonumber
\end{eqnarray}

\noindent
Since $\sum_{0\leq k\leq v_+^{-1}(n)}\rho_k\asymp n$, proceeding similarly for $\sum_{-v_-^{-1}(n)\leq k\leq l\leq0}\rho_k\rho_l\({1/\rho^2_k+1/\rho_l^2}\)$, we get~:

\begin{equation}
\label{phiplus}
\Phi_+(n)\asymp\Phi_{str}(n)+\({\sum_{-v_-^{-1}(n)\leq k\leq0\leq l\leq v_+^{-1}(n)}\rho_k\rho_l\({\sum_{s=k}^l\frac{r_s\varepsilon_s}{p_s\rho_s}}\)^2}\)^{1/2}.
\end{equation}

\noindent
In a very similar fashion :

\begin{equation}
\label{phi}\Phi(n)\asymp\Phi_{str}(n)+\({\sum_{-v_-^{-1}(n)\leq k\leq l\leq v_+^{-1}(n)}\rho_k\rho_l\({\sum_{s=k}^l\frac{r_s\varepsilon_s}{p_s\rho_s}}\)^2}\)^{1/2}.
\end{equation}

\noindent
In particular, one always has $\Phi(n)\succeq\Phi_+(n)\succeq\Phi_{str}(n)$.

\medskip
Finally an important property, called dominated variation, was proved for $\Phi^{-1}$, $\Phi_+^{-1}$ and $\Phi_{str}^{-1}$ in \cite{jb2} (lemmas 4.8 and 6.2) and will play an important role in the sequel. This condition is in fact necessitated by the structure of the main result in \cite{jb2} (theorem 2.4). Namely, for any of these function $f^{-1}$, there exists $C>0$ so that for large $x>0$ :

\begin{equation}
\label{domvar}
f^{-1}(2x)\leq Cf^{-1}(x).
\end{equation}

\noindent
As a result, for any $A>0$, for $x>0$ large enough, $f^{-1}(Ax)\asymp f^{-1}(x)$. 

\section{Proof of the results}

\noindent
- {\it Proof of proposition \ref{prop1}.} Given the remarks of the introduction, it remains to check transience when $d=2$ and $\E(\log (p_0/q_0))=0$. Following \cite{jb2}, prop. 2.5 1)$ii)$, it is enough to verify the condition on the structure function $\Phi_{str}(n)\geq\sqrt{n}(\log n)^{1/2+\varepsilon}$, $\varepsilon>0$.

\medskip
\noindent
When $\P(p_0=q_0)=1$, this is clear since $\Phi_{str}(n)\succeq n$. When $\P(p_0=q_0)<1$, introduce $S_n=\sum_{k=1}^{n}\log (q_k/p_k)$. Fixing $\varepsilon>0$, it is classical that, a.-s. for $n$ large enough~:

$$\max_{1\leq k\leq n}S_k\leq n^{1/2+\varepsilon}\mbox{ and }\min_{1\leq k\leq n}S_k\leq -n^{1/2-\varepsilon}.$$

\noindent
Hence for large $n$, $v_+(n)\leq \exp(n^{1/2+2\varepsilon})$, so $v_+^{-1}(n)\geq (\log n)^{2-9\varepsilon}$. As a result~:

$$\sum_{0\leq k\leq v_+^{-1}(n)}1/\rho_k\geq \exp(\log ^{(2-9\varepsilon)(1/2-\varepsilon)}n)\geq\exp(\log^{1-7\varepsilon}n).$$

\noindent
This gives $\Phi_{str}(n)\succeq \sqrt{n}\exp(\log^{1-\varepsilon}n)$, $\forall \varepsilon>0$, so the condition is also verified.

\bigskip
\noindent
- {\it Proof of theorem \ref{indep} i).} Let $d=1$. For the recurrence part, let $(\varepsilon_n)_{n\in\Z}$ be such that~:

$$\sum_{n\in\Z}|\varepsilon_n|/\rho_n<+\infty.$$

\noindent
This is true if $\varepsilon_n=O(\exp(-|n|^{1/2+\delta}))$, $\delta>0$, since then $\varepsilon_n/\rho_n=O(\exp(-|n|^{1/2+\delta/2}))$, which is summable on $\Z$. Next, since $\sum_{u=k}^l(r_s\varepsilon_s)/(p_s\rho_s)$ is bounded in $(k,l)$, $k\leq l$, we deduce from \eqref{phi}:

$$\Phi^2(n)\preceq \Phi^2_{str}(n)+n^2\mbox{, so }\Phi^{-1}(n)\succeq \min\{\Phi_{str}^{-1}(n),n\}.$$

\medskip
\noindent
Similary, by \eqref{phiplus}, $\Phi_+(n)\succeq\Phi_{str}(n)$, giving $\Phi_+^{-1}(n)\preceq \Phi_{str}^{-1}(n)$. In view of \cite{jb2}, theorem 2.4, recurrence is equivalent to~:

$$\sum_{n\geq1}\frac{1}{n^2}\frac{(\Phi^{-1}(n))^2}{\Phi_+^{-1}(n)}=+\infty.$$

\noindent
It is thus sufficient to show the divergence of~:

$$\sum_{n\geq1}\frac{1}{n^2}\frac{(\min\{\Phi_{str}^{-1}(n),n\})^2}{\Phi_{str}^{-1}(n)}=\sum_{n\geq1}\frac{1}{n}\min\left\{{\frac{\Phi_{str}^{-1}(n)}{n},\frac{n}{\Phi_{str}^{-1}(n)}}\right\}.$$

\noindent
Since $\Phi_{str}^{-1}(n)$ verifies dominated variation \eqref{domvar}, it is equivalent to check that~:

\begin{equation}
\label{dvser}
\sum_{n\geq1}\min\left\{{\frac{\Phi_{str}^{-1}(K^n)}{K^n},\frac{K^n}{\Phi_{str}^{-1}(K^n)}}\right\}=+\infty.
\end{equation}

\noindent
We shall in fact prove that the general term does not go to $0$. 

\begin{lemme} 

$ $

\noindent
Introduce for $n\geq0$~:

$$w_+(n)=\sum_{0\leq k\leq n}1/\rho_k\mbox{ and }w_-(n)=(1/a_0)\sum_{-n-1\leq k\leq -1}1/\rho_k.$$

$ $

\noindent
A.-s., $\displaystyle\limsup \frac{w_+(n)}{v_+(n)}=+\infty$ and $\displaystyle\limsup\min\left\{{\frac{v_+(n)}{w_+(n)},\frac{v_-(n)}{v_+(n)},\frac{w_+(n)}{w_-(n)}}\right\}=+\infty$.\end{lemme}

\noindent
\textit{Proof of the lemma~:}

\noindent
Let us prove the second point, the first one being easier. Fix a sequence $(k_n)$ with $k_n/k_{n-1}^2\rightarrow+\infty$. For $n\in\Z$, let $S_{n}=\log \rho_n$. Let~:

$$U_n=\max_{k\in[-k_n,-k_{n-1})}(S_k-S_{-k_{n-1}}),~V_n=\max_{k\in(k_{n-1},k_{n}]}(S_k-S_{k_{n-1}}).$$

\noindent
In the same way, introduce~:

$$W_n=\min_{k\in(k_{n-1},k_{n}]}(S_k-S_{k_{n-1}}),~X_n=\min_{k\in[-k_{n},-k_{n-1})}(S_k-S_{-k_{n-1}}).$$

\noindent
Set $\sigma^2=\mbox{Var}(\log (q_0/p_0))>0$ and let $c=2+2|\log(1/\eta)|$. Using functional convergence to standard Brownian motion $(B_t)_{t\in[-1,1]}$, we have, as $n\rightarrow+\infty$~:

\begin{eqnarray}
&&\P\({\frac{U_n}{\sigma\sqrt{k_n-k_{n-1}}}\geq 1+\frac{V_n}{\sigma\sqrt{k_n-k_{n-1}}}\geq2-\frac{W_n}{\sigma\sqrt{k_n-k_{n-1}}}\geq3-\frac{X_n}{\sigma\sqrt{k_n-k_{n-1}}}}\)\nonumber\\
&\longrightarrow& \P\({\max_{t\in[-1,0]}B_t\geq1+\max_{t\in[0,1]}B_t\geq2-\min_{t\in[0,1]}B_t\geq3-\min_{t\in[-1,0]}B_t}\)=:\alpha>0.\nonumber
\end{eqnarray}

\noindent
Using independence and Borel-Cantelli 2, almost-surely the event appearing in the first probability is realized for infinitely many $n$. For such a $n$, we have~:

\begin{eqnarray}
v_+(k_n)&\geq& \exp(-ck_{n-1}+\max_{k\in(k_{n-1},k_n]}(S_k-S_{k_{n-1}}))\nonumber\\
&\geq&\exp(-ck_{n-1}+\sigma\sqrt{k_n-k_{n-1}}-\min_{k\in(k_{n-1},k_n]}(S_k-S_{k_{n-1}}))\nonumber\\
&\geq&\exp(-3ck_{n-1}+\sigma\sqrt{k_n-k_{n-1}}-\min_{k\in[1,k_n]}S_k)\nonumber\\
&\geq&\frac{w_+(k_n)}{k_n}\exp(-3ck_{n-1}+\sigma\sqrt{k_n-k_{n-1}}).\nonumber\end{eqnarray}

\noindent
If $n$ is large, $v_+(k_n)/w_+(k_n)\geq \exp(\sigma\sqrt{k_n}/2)/k_n$, which is arbitrary large. Similar lower bounds are proved for $v_-(k_n)/v_+(k_n)$ and $w_+(k_n)/w_-(k_n)$.

\fin

\bigskip
\noindent
We conclude the proof of the recurrence part. Using the first point of the lemma, a.-s., for any $A>0$, infinitely often, $w_+(n)\geq Av_+(n)$, i.e. $w_+(v_+^{-1}(v_+(n)))\geq A v_+(n)$. Taking $l$ so that $l\leq v_+(n)<l+1$, we obtain~:

$$\Phi_{str}(l+1)\geq\sqrt{l+1}\sqrt{Al}\geq l\sqrt{A}.$$

\medskip
\noindent
Hence, $\Phi_{str}^{-1}(l\sqrt{A})\leq l+1$. This implies that, a.-s., $\liminf\Phi_{str}^{-1}(n)/n=0$. Using now the second point of the lemma, a.-s., for any $A>1$, i.o., both $v_-(n)\geq v_+(n)\geq Aw_+(n)$ and $w_+(n)\geq w_-(n)$. For such a $n$, first $v_-^{-1}(v_+(n))\leq n$. Then~:

\begin{eqnarray}
\frac{v_+(n)}{\sum_{-v^{-1}(v_+(n))\leq k\leq v_+^{-1}(v_+(n))}1/\rho_k}&\geq&\frac{v_+(n)}{\sum_{-n\leq k\leq n}1/\rho_k}\geq \frac{v_+(n)}{a_0w_-(n)+w_+(n)}\nonumber\\
&\geq&\frac{v_+(n)}{w_+(n)}(1+a_0)^{-1}\geq\frac{A}{1+a_0}.\nonumber
\end{eqnarray}

\noindent
Choosing $l$ such that $l\leq v_+(n)<l+1$, we obtain~:

$$\frac{l+1}{\Phi_{str}(l)}\geq\frac{A}{1+a_0}.$$

\noindent
Hence a.-s., $\limsup n/\Phi_{str}(n)=+\infty$ and so, in the same way as before, $\limsup \Phi_{str}^{-1}(n)/n=+\infty$. Let now $b_n=\Phi_{str}^{-1}(K^n)/K^n$. Using dominated variation \eqref{domvar}, a.-s., the previous results furnish that $\liminf b_n=0$, $\limsup b_n=+\infty$. Dominated variation also implies that, a.-s. for a constant $H>1$, for all $n$, we have $b_n/b_{n+1}\in[1/H,H]$. Thus $b_n\in[1,H]$ for infinitely many $n$. For such a $n$, $\min\{b_n,1/b_n\}=1/b_n\geq 1/H$ and this shows \eqref{dvser}.

\bigskip
Turn now to the proof of the transience result of item $i)$. By \cite{jb2}, prop. 2.5 $1)$, it is enough to prove that $\sum_{n\geq1}1/\Phi(n)<+\infty$. We shall assume that the $\varepsilon_n$ are non-negative for large $n\in\Z$, with $\sum_{n\in\Z}\varepsilon_n/\rho_n=+\infty$ and~:

\begin{equation}
\label{reste}\sum_{n\geq N_0}\({n\sum_{-v_-^{-1}(n)\leq k\leq v_+^{-1}(n)}\varepsilon_k/\rho_k}\)^{-1}<+\infty,
\end{equation}

\noindent
where $N_0$ is chosen large enough so that the denominator of the generic term above is $>0$ for $n\geq N_0$. This is verified if there exists $\delta>0$ so that $\varepsilon_n\geq \exp(-|n|^{1/2-\delta})$, for large $n\in\Z$. Indeed, as in the proof of prop. \ref{prop1}, a.-s., for any $\varepsilon>0$, for $n$ large enough~:

$$\min_{0\leq k\leq v_+^{-1}(n)}\rho_k\leq \exp(-\log^{1-\varepsilon}n).$$

\noindent
If $u_n$ is the first point in $[0,v_+^{-1}(n)]$ realizing the minimum of $\rho_k$ on this interval, then $u_n\rightarrow+\infty$ and, fixing $\varepsilon>0$, $u_n\leq (\log n)^{2+\varepsilon}$, for large $n$. Thus for $n$ large enough~:

$$\varepsilon_{u_n}/\rho_{u_n}\geq \exp(-u_n^{1/2-\delta})\exp(\log^{1-\varepsilon}n)\geq \exp(\log^{1-2\varepsilon}n),$$

\noindent
if $\varepsilon>0$ was chosen small enough. Then $\forall \varepsilon>0$, for large $n$~:

$$\sum_{-v_-^{-1}(n)\leq k\leq v_+^{-1}(n)}\varepsilon_k/\rho_k\geq \frac{1}{2}\varepsilon_{u_n}/\rho_{u_n}\geq \exp(\log^{1-2\varepsilon}n).$$

\noindent
This shows \eqref{reste}. Notice that because of monotonicity for large $n$, condition \eqref{reste} is equivalent to~:

\begin{equation}
\label{condens}
\sum_{n\geq N_0}\({\sum_{-v_-^{-1}(K^n)\leq k\leq v_+^{-1}(K^n)}\varepsilon_k/\rho_k}\)^{-1}<+\infty.
\end{equation}

\noindent
Now, because $\Phi$ is increasing, we have to check that $\sum_{n\geq1}K^n/\Phi(K^n)<+\infty$. Observe that for large $n$~:

\begin{eqnarray}
\Phi^2(K^n)&\geq& \sum_{\overset{-v_-^{-1}(K^n)\leq k<-v_-^{-1}(K^{n-1})}{v_+^{-1}(K^{n-1})<l\leq v_+^{-1}(K^n)}}\rho_k\rho_l\({\sum_{-v_-^{-1}(K^{n-1})\leq s\leq v_+^{-1}(K^{n-1})}\frac{r_s\varepsilon_s}{p_s\rho_s}}\)^2\nonumber\\
&\succeq&\({\sum_{-v_-^{-1}(K^{n-1})\leq s\leq v_+^{-1}(K^{n-1})}\frac{\varepsilon_s}{\rho_s}}\)^2 \sum_{\overset{-v_-^{-1}(K^n)\leq k<-v_-^{-1}(K^{n-1})}{v_+^{-1}(K^{n-1})<l\leq v_+^{-1}(K^n)}}\rho_k\rho_l\nonumber\\
&\succeq&K^{2n}\({\sum_{-v_-^{-1}(K^{n-1})\leq s\leq v_+^{-1}(K^{n-1})}\frac{\varepsilon_s}{\rho_s}}\)^2,\nonumber\end{eqnarray}

\noindent
since~:

\begin{eqnarray}K^n&\geq&\sum_{k=v_+^{-1}(K^{n-1})+1}^{v_+^{-1}(K^n)}\rho_k=v_+(v_+^{-1}(K^n))-v_+(v_+^{-1}(K^{n-1}))\geq\frac{K^n}{1+1/\eta}-K^{n-1}\geq\frac{K^n}{2(1+1/\eta)},\nonumber
\end{eqnarray} 

\noindent
with similar inequalities for $\sum_{-v_-^{-1}(K^n)\leq k<-v_-^{-1}(K^{n-1})}\rho_k$. Hence $\sum_{n\geq1}K^n/\Phi(K^n)<+\infty$, by \eqref{condens}.

\bigskip
\noindent
- {\it Proof of theorem \ref{indep} ii).} We now assume that the $(\varepsilon_n)_{\in\Z}$ are independent random variables, independent from the $(p_n/q_n)_{n\in\Z}$, such that there exists $\delta>0$ such that for large $n\in\N$ and all $y\in\R$, $\P(\varepsilon_n\not\in B(y,\exp(-n^{1/2-\delta})))\geq n^{-\delta/5}$. 

\medskip
\noindent
First of all, as above~:

$$\sum_{v_+^{-1}(K^{n-1})<l\leq v_+^{-1}(K^n)}\rho_l\asymp K^n\mbox{ and }\sum_{-v_-^{-1}(K^n)\leq l\leq-1}\rho_l\asymp K^n.$$

\noindent
Introducing then the probability measure on $\Z^2$~:

$$\nu_n=\frac{1}{Z_n}\sum_{-v_-^{-1}(K^n)\leq k\leq -1,v_+^{-1}(K^{n-1})< l\leq v_+^{-1}(K^n)}\rho_k\rho_l\delta_{(k,l)},$$

\noindent
the normalizing constant $Z_n$ verifies $Z_n\asymp K^{2n}$.

\medskip
Let $\delta/2>\gamma>2\delta/5$. Next choose $\varepsilon>0$ so that $\gamma+(2+\varepsilon)(1/2-\delta)\leq 1-3\varepsilon$ and $(2+\varepsilon)\delta/5<\gamma<\delta/2$. Let $u_n$ in $[0,v_+^{-1}(K^{n-1})]$ be the first minimum of $\rho_k$ on this interval. A.-s. for large $n$, we have $n^{2+\varepsilon}\geq v_+^{-1}(K^n)\geq u_n\geq n^{2-\varepsilon}$ and~:

$$\frac{1}{\rho_{u_n}}\geq \exp(n^{1-\varepsilon}).$$

\noindent
The support of the law of $\varepsilon_n$ is included in $[-1/\eta,1/\eta]$. Let $C_0=2/\eta+1$ and $l(n)=n^{\gamma}$. Starting from $u_{n,0}=u_n$, choose recursively in a decreasing order in the interval $[0,u_n]$, points $u_{n,0}>u_{n,1}>\cdots>u_{n,l(n)}>0$, so that $u_{n,l+1}$ is the closest point on the left side of $u_{n,l}$ with~:

$$\frac{1}{4}\frac{C_0}{\rho_{u_{n,l}}}\geq\frac{1}{4}\frac{\exp(-u_{n,l}^{1/2-\delta})}{\rho_{u_{n,l}}}\geq \frac{C_0}{\rho_{u_{n,l+1}}} \geq\frac{1}{4K}\frac{\exp(-u_{n,l}^{1/2-\delta})}{\rho_{u_{n,l}}},~0\leq l<l(n).$$

\medskip
\noindent
Recall for this that $\rho_k/\rho_{k+1}\in[1/K,K]$ and notice that $\rho_{u_{n,l}}$ is increasing, $0\leq l\leq l(n)$, and indeed well defined and staying inside the interval $[0,u_{n,0}]$ since uniformly in $0\leq l< l(n)$~:

\begin{eqnarray}
\label{1mep}
\rho_{u_{n,l+1}}\leq (4C_0K) \rho_{u_{n,l}}e^{u_{n,l}^{1/2-\delta}}&\leq& (4C_0K)^{l(n)} \rho_{u_{n,0}}\exp(l(n)u_{n,0}^{1/2-\delta})\nonumber\\
&\leq& \exp(n^{\gamma}\log(4C_0K)-n^{1-\varepsilon}+n^{\gamma}n^{(2+\varepsilon)(1/2-\delta)})\nonumber\\
&\leq&\exp(-n^{1-2\varepsilon}),
\end{eqnarray}

\noindent for $n$ large enough. In particular, $u_{n,l(n)}\rightarrow+\infty$.

\medskip
We now reason conditionally to the $(p_n/q_n)_{n\in\Z}$ and make a measurable construction. The probability $\P$ hence only refers to the $(\varepsilon_n)_{n\in\Z}$. Consider the event $A_n$, where $(\varepsilon_i)_{i\in\Z}\in A_n$ if ~: 

\begin{equation}
\label{nun}
\nu_n\left\{{(k,l)\in\Z^2,~\left|{\sum_{k\leq u\leq l}\frac{\varepsilon_u}{\rho_u}}\right|\leq \frac{C_0}{\rho_{u_{n,l(n)}}}}\right\}\geq 3/4.
\end{equation}

\medskip
\noindent
Let us write~:

$$\P(A_n)=\E\({\E(1_{(\varepsilon_k)_{k\in\Z}\in A_n}~|~(\varepsilon_l)_{l\not=u_{n,s},0\leq s<l(n)})}\).$$

\medskip
\noindent
Fix $(\varepsilon_l)_{l\not=u_{n,s},0\leq s<l(n)}$ and suppose that there exists $(\varepsilon_{u_{n,s}})_{0\leq s<l(n)}$ so that $\varepsilon=(\varepsilon_i)_{i\in\Z}\in A_n$. Let~:

$$(\varepsilon'_{u_{n,s}})_{0\leq s<l(n)}\not\in B\({\varepsilon_{u_{n,l(n)-1}},\exp(-u_{n,l(n)-1}^{1/2-\delta})}\)\times\cdots\times B\({\varepsilon_{u_{n,0}},\exp(-u_{n,0}^{1/2-\delta})}\).$$

\medskip
\noindent
We call $\varepsilon'$ the point obtained from $\varepsilon$ by replacing $\varepsilon_{u_{n,l}}$ by $\varepsilon'_{u_{n,l}}$, $0\leq l<l(n)$, without changing the other coordinates. We claim that $\varepsilon'\not\in A_n$. Indeed, let $0\leq p< l(n)$ be the smallest index so that~:

$$\varepsilon'_{u_{n,p}}\not \in B\({\varepsilon_{u_{n,p}},\exp(-u_{n,p}^{1/2-\delta})}\).$$

\noindent
Take $(k,l)$ in the set involved in \eqref{nun}. Notice that $k<0$ and $v_+^{-1}(K^{n-1})<l$. Using the hypothesis on the laws of the $\varepsilon_n$~:

\begin{eqnarray}
|\sum_{k\leq u\leq l}\frac{\varepsilon'_u}{\rho_u}|&=&|\frac{\varepsilon'_{u_{n,p}}-\varepsilon_{u_{n,p}}}{\rho_{u_{n,p}}}+\sum_{p<q< l(n)}\frac{\varepsilon'_{u_{n,q}}-\varepsilon_{u_{n,q}}}{\rho_{u_{n,q}}}+\sum_{k\leq u\leq l}\frac{\varepsilon_u}{\rho_u}|\nonumber\\
&\geq&\frac{\exp(-u_{n,p}^{1/2-\delta})}{\rho_{u_{n,p}}}-\sum_{p<q< l(n)}\frac{C_0}{\rho_{u_{n,q}}}-\frac{C_0}{\rho_{u_{n,l(n)}}}\nonumber\\
&\geq &4\frac{C_0}{\rho_{u_{n,p+1}}}-\frac{C_0}{\rho_{u_{n,p+1}}}\sum_{r\geq0}4^{-r}\geq 2\frac{C_0}{\rho_{u_{n,p+1}}}\geq2\frac{C_0}{\rho_{u_{n,l(n)}}}.
\end{eqnarray}

\medskip
\noindent
This furnishes~:

\begin{equation}
\label{nun2}
\nu_n\left\{{(k,l)\in\Z^2,~\left|{\sum_{k\leq u\leq l}\frac{\varepsilon'_u}{\rho_u}}\right|\leq \frac{C_0}{\rho_{u_{n,l(n)}}}}\right\}\leq 1/4,
\end{equation}

\medskip
\noindent
hence proving that $\varepsilon'\not\in A_n$. As a result~:

\begin{eqnarray}
\E(1_{(\varepsilon_k)_{k\in\Z}\in A_n}~|~(\varepsilon_l)_{l\not=u_{n,s},0\leq s<l(n)})&\leq&(1-u^{-\delta/5}_{n,l(n)-1})\cdots(1-u^{-\delta/5}_{n,0})\nonumber\\
&\leq&(1-u_{n}^{-\delta/5})^{l(n)}\leq \exp(-n^{\gamma-(2+\varepsilon)\delta/5}),\nonumber
\end{eqnarray}

\noindent
for large $n$. This gives $\P(A_n)\leq \exp(-n^{\gamma-(2+\varepsilon)\delta/5})$, so $\sum \P(A_n)<\infty$, since $\gamma-(2+\varepsilon)\delta/5>0$. Making use of Borel-Cantelli 1, for almost-all $(\varepsilon_i)_{i\in\Z}$, for $n$ large enough~:

\begin{equation}
\label{nun3}
\nu_n\left\{{(k,l)\in\Z^2,~\left|{\sum_{k\leq u\leq l}\frac{\varepsilon_u}{\rho_u}}\right|\geq \frac{C_0}{\rho_{u_{n,l(n)}}}}\right\}\geq 1/4.
\end{equation}

\medskip
\noindent
We now conclude. Write ${\cal E}_n(f(X,Y))$ for $\sum_{(k,l)\in\Z^2}f(x,y)d\nu_n(k,l)$. Using \cite{jb2}, prop. 2.5, we prove that $\sum_{n\geq1}1/\Phi(n)<+\infty$, or equivalently $\sum_{n\geq1}K^n/\Phi(K^n)<+\infty$, as $\Phi$ is increasing. We have~:

\begin{eqnarray}
\sum_{n\geq1}\frac{K^n}{\Phi(K^n)}&\leq&\sum_{n\geq1} K^n\left({\sum_{-v_-^{-1}(K^n)\leq k\leq l\leq v_+^{-1}(K^n)}\rho_k\rho_l\({\sum_{k\leq u\leq l}\frac{\varepsilon_u}{\rho_u}}\)^2}\)^{-1/2}\nonumber\\
&\preceq&\sum_{n\geq1} \({{\cal E}_n\[{\({\sum_{X\leq u\leq Y}\frac{\varepsilon_u}{\rho_u}}\)^2}\]}\)^{-1/2}.\nonumber
\end{eqnarray}

\medskip
\noindent
From the previous results, a.-s., there exists $N_0$ so that for $n\geq N_0$~:

$${\cal E}_n\[{\({\sum_{X\leq u\leq Y}\frac{\varepsilon_u}{\rho_u}}\)^2}\]\geq\frac{1}{4}\frac{C_0^2}{\rho^{2}_{u_{n,l(n)}}}.$$

\medskip
\noindent
From \eqref{1mep}, $\sum_{n\geq1}\rho^{-1}_{u_{n,l(n)}}\preceq\sum_{n\geq1}\exp(-n^{1-2\varepsilon})<+\infty$. This shows $\sum_{n\geq1}K^n/\Phi(K^n)<+\infty$ and completes the proof of the theorem. \fin

\section{Concluding remarks}
This is an informal section. Let us fix the context of theorem \ref{indep}. As in \cite{jb1,jb2}, introduce the random times $0=\sigma_0<\tau_0<\sigma_1<\tau_1<\cdots$, where $\tau_k=\min\{n>\sigma_k~|~S_n\not\in\Z\times\{0\}\}$ and $\sigma_{k+1}=\{n>\tau_k~|~S_n\in\Z\times\{0\}\}$. Set $D_n=S_{\sigma_n}-S_{\sigma_{n-1}}$ and recall that the stratification of the environment implies that the $(D_n)_{n\geq 1}$ are independent and identically distributed.

\medskip
\noindent
$i)$  The condition $\sum_{n\in\Z}|\varepsilon_n|/\rho_n<+\infty$, implying recurrence in theorem \ref{indep} $i)$, looks like the condition $\sum_{n\in\Z}1/\rho_n<+\infty$ defining the half-pipe of \cite{jb2}, section 7.3, but is in fact of different nature. This last condition, not true here, implied that $D_1$ is integrable and that the random walk in this case is transient if and only if~:

\begin{equation}
\label{halfpipe}
\sum_{s\in\Z}\frac{r_s\varepsilon_s}{p_s\rho_s}\not=0.
\end{equation}

\noindent
Here, under the condition $\sum_{n\in\Z}|\varepsilon_n|/\rho_n<+\infty$, whatever the value of $\sum_{s\in\Z}r_s\varepsilon_s/(p_s\rho_s)$, the random walk is recurrent.

\medskip
\noindent
$ii)$ We do not expect that a complete characterization of transience in such a simple form as \eqref{halfpipe} is possible in general. When $\varepsilon_n>0$, $\forall n\in\Z$, setting~:

$$C(n)=\sum_{-v_-^{-1}(n)\leq k\leq v_+^{-1}(n)}\varepsilon_k/\rho_k,$$

\noindent
mention without proof that the following condition is equivalent to recurrence~:

$$\sum_{n\geq1}\frac{1}{n^2}\min\{\Phi_{str}^{-1}(n),(C^{-1}(n))^2/\Phi_{str}^{-1}(n)\}=+\infty.$$

\noindent
Notice that the formulation only has a meaning when the $\varepsilon_n$ are positive (because of $C^{-1}(n)$). 

\medskip
\noindent
$iii)$ Standardly, see Gnedenko-Kolmogorov \cite{gk}, the correct normalization for $D_1+\cdots+D_n$ is related to the behaviour at the origin of the characteristic function of $D_1$. Here, as a consequence of \cite{jb2}, almost-surely, the $(D_1+\cdots+D_n)_{1\leq n\leq N}$ have a standard empirical deviation $d_N$ (the square root of the empirical variance) around its empirical mean $m_N$ given by $\Phi_+(N)$ and a global size given by $\Phi(N)$ (i.e. $|m_N|+d_N$). Taking $\varepsilon_n=1$, one may deduce that the correct normalization for $\sigma_n$ has order $\tilde{\Phi}(n)$, where~:

$$\tilde{\Phi}(n)\asymp\({\sum_{-v_-^{-1}(n)\leq k\leq l\leq v_+^{-1}(n)}\rho_k\rho_l\({\sum_{s=k}^l\frac{1}{\rho_s}}\)^2}\)^{1/2}.$$

\noindent
From this it is not difficult to infer that $\tilde{\Phi}^{-1}(n)\asymp \Psi^{-1}(n)$, where $\Psi(n)=n\sum_{-v_-^{-1}(n)\leq k\leq v_+^{-1}(n)}1/\rho_k$.

\medskip
In the present context, the correct normalization for $S_n$ should then be $(\Phi\circ{\Psi}^{-1}(n),(\log n)^2)$. Of course some work is required to make it rigourous, for example to prove tightness or some limit theorem in law (making some stationarity hypothesis on the $(\varepsilon_n)_{n\in\Z}$). Recall that in the vertically flat case, when the $(\varepsilon_n)_{n\in\Z}$ are $i.i.d.$ centered and non-constant, then $\Phi(n)\asymp n^{3/2}$ and ${\Psi}(n)\asymp n^2$, hence recovering the size $n^{3/4}$ for the first coordinate, as shown in \cite{DP}.

\bigskip
\noindent
{\bf{Acknowledgments.}} We would like to thank Jacques Printems and Alexis Devulder for useful discussions. We also thank Fran\c coise P\`ene for bringing reference \cite{koch} to our attention.

\providecommand{\bysame}{\leavevmode\hbox to3em{\hrulefill}\thinspace}

\bigskip
{\small{\sc{Laboratoire d'Analyse et de Math\'ematiques Appliqu\'ees,
  Universit\'e Paris-Est, Facult\'e des Sciences et Technologies, 61, avenue du G\'en\'eral de Gaulle, 94010 Cr\'eteil Cedex, FRANCE}}}

\it{E-mail address~:} {\sf julien.bremont@u-pec.fr}

\end{document}